\documentclass[review,10pt]{elsarticle}
\usepackage[numbers]{natbib}
 \usepackage{mathrsfs}
\usepackage{graphicx}
\usepackage{epstopdf}
\usepackage{amssymb}
\usepackage{color}
\usepackage{amsmath}
\usepackage{float}
\makeatletter\@addtoreset{equation}{section}

\newtheorem{thm}{Theorem}[section]

\newtheorem{lem}[thm]{Lemma}

\newtheorem{rem}[thm]{Remark}

\numberwithin{equation}{section}
\numberwithin{table}{section}
\numberwithin{figure}{section}

\makeatletter\@addtoreset{equation}{section}

\newenvironment{proofed}[1]{\par \textbf{Proof}\quad #1}{\hfill \textbf{} $\Box$ }

\newtheorem{example}{\bf{Example}}[section]

\numberwithin{equation}{section}
\numberwithin{table}{section}
\numberwithin{figure}{section}

\begin{document}


\begin{frontmatter}
\title{A second-order difference scheme for the  time fractional substantial diffusion equation}

\author[seu]{Zhaopeng Hao}
\ead{230139121@seu.edu.cn}
\author[seu]{Wanrong Cao\corref{cor1}}
\ead{wrcao@seu.edu.cn}
\cortext[cor1]{Corresponding author. Tel/Fax: +86 25 52090590}
\author[purdue]{Guang Lin}
\ead{guanglin@purdue.edu}

\address[seu]{Department of Mathematics, Southeast University, Nanjing 210096, P.R.China.}
\address[purdue]{Department of Mathematics, School of Mechanical Engineering, Purdue University, West Lafayette, IN 47907, USA.}

\begin{abstract}
In this work,  a second-order approximation of  the fractional substantial derivative is presented  by considering a modified shifted substantial Gr\"{u}nwald formula and its asymptotic expansion. Moreover,  the proposed approximation is applied  to a fractional diffusion equation with fractional substantial derivative in time. With the use of the  fourth-order compact scheme in space, we give a fully discrete Gr\"{u}nwald-Letnikov-formula-based compact difference scheme and prove its stability and convergence by the energy method under smooth assumptions. In addition, the problem with nonsmooth solution is also discussed,  and an improved algorithm is proposed to deal with the singularity of the fractional substantial derivative.     Numerical examples show the reliability and efficiency of the scheme.
\end{abstract}


\end{frontmatter}

\section{Introduction}
Anomalous sub-diffusion process,  commonly described by the continuous time random walks (CTRWs), and also known as non-Brownian sub-diffusion, arises in numerous physical,chemical and biological systems; see \cite{BronsteinIK2009,Cabal2006-nature-anomulas-diffusion,CometeDT2005,KleinhansF2007}.  The Feynman-Kac formula named after Richard Feynman and Mark Kac, establishes a link between parabolic partial differential equations (PDEs) and Brownian functionals. To figure out the probability density function (PDF) of some non-Brownian functionals, the fractional Feynman-Kac equation has been derived in \cite{CarmiB2011,CarmiTB2010,FriedrichJBE2006,SokolovM2003,TurgemanCB2009}. The non-Brownian functionals can be defined by $A=\int_0^tU(x(\tau))\mathrm{d}\tau$, where $x(t)$ is the trajectory of a non-Brownian particle and different choices of $U(x)$ can depict diverse systems. In particular, if taking $U(x)\equiv 0$, the fractional Feynman-Kac equation reduces to the well-known fractional Fokker-Planck equation; see \cite{CarmiB2011,CarmiTB2010,MetzlerK2000} for details.  {\color{red} L$\acute{e}$vy walks  give a proper
dynamical description in the superdiffusive domain, where the
temporal and spatial variables of L$\acute{e}$vy walks are strongly correlated and the PDFs of waiting
time and jump length are spatiotemporal coupling \cite{SokolovM2003}. Thus, an important operator, fractional substantial derivative has been proposed to describe 
 the CTRW models with coupling PDFs.  This spatiotemporal coupling operator was also  presented in  \cite{FriedrichJBE2006}, where  the CTRW model with position-velocity coupling PDF was discussed.   Recently,  Carmi and Barkai \cite{CarmiB2011} also used   the substantial derivative to derive the forward and backward fractional Feynman-Kac
equations. Due to its potential properties and wide application, the fractional substantial derivative has attracted many scholars' attention; see \cite{CaoLC2015,ChenD2014} and references therein.}

The fractional substantial derivative operator of order $\alpha\;(n-1<\alpha<n)$ is defined by \cite{ChenD2014}
\begin{equation}\label{eq:def-sub}
 {_{a}}D_t^{\alpha,\lambda}f(t)={_{a}}D_t^{n,\lambda} [{_{a}}I_t^{n-\alpha,\lambda}f](t), \quad {_{a}}I_t^{\alpha,\lambda}f(t)=\frac{1}{\Gamma(\alpha) }\int_a^t\frac{ f(\xi)}{(t-\xi)^{1-\alpha} e^{\lambda(t-\xi) }} d\xi;
\end{equation}
where  $\alpha>0$, $\lambda$ is a constant or a function independent of variable $t$, ${_{a}}I_t^{\alpha,\lambda}$ denotes the fractional substantial operator, and
$${_{a}}D_t^{n,\lambda}=(\frac{d}{d t }+\lambda)^{n}=(\frac{d}{d t }+\lambda)(\frac{d}{d t }+\lambda)\cdots(\frac{d}{d t }+\lambda). $$
{\color{red}It is noted that if $\lambda$ is a nonnegative constant, then the fractional substantial derivative is equivalent to the Riemnan-Liouville tempered derivative defined in \cite{BaeumerMeer2010,Buschman1972,Cartea2007}, and taking $\lambda=0$ in \eqref{eq:def-sub} leads to the left Riemnan-Liouville derivative. Meanwhile, to obtain  the non-Brownian functionals, whose path integrals {\color{red}are} over L$\acute{e}$vy trajectories, the {\color{red}space-fractional} Fokker-Planck equation and the tempered space fractional diffusion equations have  been widely used; see \cite{Hanert2014-tempered,Kullberg2012-tempered,MeerschaertZB2008}.}

The current work is devoted to proposing a second-order Gr\"{u}nwald-Letnikov-formula-based approximation for the fractional substantial derivative \eqref{eq:def-sub}, and applying it to derive a high-order fully discrete scheme for the time fractional substantial diffusion equation (TFSDE)
\begin{eqnarray}
&& {\color{red}{_0}D_t^{\alpha,\lambda} }u (\mathbf{x},t)=\Delta u(\mathbf{x},t)+F(\mathbf{x},t),  \quad \mathbf{x}\in \Omega,\quad  t\in (0, T], \label{c1}\\
&&u(\mathbf{x},0)=u_0(\mathbf{x}), \quad \mathbf{x} \in {\Omega},\label{cc0}\\
 &&  u(\mathbf{x},t)|_{\mathbf{x}\in \partial{\Omega}}=\phi(\mathbf{x},t),  \quad t\in (0, T],\label{cc1}
  \end{eqnarray}
where  $\Delta $  is the Laplacian operator, $\mathbf{x}$ denotes the one-dimensional or two-dimentional space variable, $\partial{\Omega}$ is the boundary of domain $\Omega$, $F (\mathbf{x},t),$ $u_0(\mathbf{x})$ and $\phi(\mathbf{x},t)$ are given functions
; ${\color{red}{_0}D_t^{\alpha,\lambda} }$ is the substantial derivative defined by \eqref{eq:def-sub}, and $0< \alpha \leq 1.$


The main novelty of our paper is the derivation of a second-order operator, which is based on the modified definition of the Gr\"{u}nwald derivative (see section 3.4, \cite{OldhamS1974}), for the approximation of the fractional substantial derivative. The modified Gr\"{u}nwald derivative is defined by \cite{OldhamS1974}
\begin{eqnarray}\label{GL-M}
&& {}_0^{GL}D_t^\alpha f(t)\equiv\lim_{\tau \rightarrow 0}\frac 1{\tau^\alpha}\sum_{k=0}^{[t/\tau]}g_k^\alpha f(t-k\tau+\frac{\alpha}{2}\tau), \quad g_{k}^{\alpha}=(-1)^k\left(
\begin{array}{l}
\alpha\\
k
\end{array}\right).\qquad
\end{eqnarray}
Actually, if dropping the term $\displaystyle \frac{\alpha}{2}\tau$ off in the right side of above definition, one gets the original definition of the Gr\"{u}nwald derivative. The advantage of the modified Gr\"{u}nwald derivative is that it permits the design of more efficient algorithms to approximate the Riemann-Liouville fractional derivatives than using the shifted Gr\"{u}nwald-Letnikov formula directly.  In this paper, by developing the modified Gr\"{u}nwald derivative and the shifted Gr\"{u}nwald-Letnikov formula to the fractional substantial derivative,  a modified shifted substantial Gr\"{u}nwald formula and  its asymptotic expansion  are presented. Based upon the asymptotic expansion, a second-order approximation of  the fractional substantial derivative is derived.

There are also some other approaches for the approximation of fractional derivatives, such as the $L1$ approximation \cite{LinXu07,SunWu06}, the fractional linear multi-step methods (FLMMs) developed by Lubich \cite{ZengLLTI2015}, the $L2$ approximation with using superior convergence \cite{Alikhanov2015,GaoSun2015-superconvergence}, etc.   However, {\color{red}to the  best of authors' knowledge}, very limited  work has been presented for the fractional substantial derivative.  Chen and Deng \cite{ChenD2014} extended the {\color{red}$p$}-th order FLMMs \cite{Lubich1986} to approximate the fractional substantial derivative, and applied it to solve the fractional  Feyman-Fac equation \cite{DengCB2015-JSC-substantial} lately. Very recently, Chen and Deng \cite{ChenDeng2014-arxiv-substantial-tempered} proposed some algorithms for the equation with  the fractional substantial derivative in time  and  the tempered fractional derivatives in space, in which the numerical stability and error estimate have been given for a scheme with {\color{red}the} first-order accuracy in time and {\color{red}the} second-order accuracy in space.

The main goal of our paper is to construct a second-order approximation for the time fractional substantial derivative, and subsequently to solve the TFSDE \eqref{c1}-\eqref{cc1} by combining the existed fourth-order  {\color{red}compact approximation} for the  space  derivatives \cite{sun2005}, and establish  the numerical stability and error estimate of the derived fully discretized scheme.

In this work, we assume that the solution to the underlying equation
satisfies suitable regularity requirements.   The assumption can be satisfied in certain conditions,   while
it may not hold for  many time-fractional differential equations; see related  discussion  for the case $\lambda=0$ in \cite{CaoZG2015,JinLZ13,ZengLLT13}. To circumvent the requirement of high regularity of the solution, we apply starting quadrature to add correction terms in the proposed scheme. The starting quadrature was first developed in \cite{Lubich1986}, and {\color{red} has} been used to deal with problems with nonsmooth solution; see \cite{DengCB2015-JSC-substantial,ZengLLTI2015}. The validity of the proposed algorithm is illustrated in Example 6.3 by solving a two-dimensional time fractional substantial diffusion equation.

The remainder of this paper is organized as follows. In Section 2,  a second-order operator for the approximation of the fractional substantial derivative is derived.   In Section 3,  the proposed approximation is applied to TFSDE \eqref{c1}-\eqref{cc1}, and   a fully discretized scheme is derived  by combining the fourth-order compact formula in space. In Section 4, we give a discrete prior estimate for the numerical solution, and then prove the convergence and stability of the proposed scheme. The behavior of our proposed scheme when applied to solve problems with non-smooth solution is further  discussed and the improved algorithm is proposed   in Section 5.  Numerical results are presented in Section 6.  Some concluding remarks are included in the final section.

\section{A second-order approximation for the fractional substantial  derivative}
\setcounter{equation}{0}
In this section, we first present a second-order approximation for the fractional substantial derivative.
For $0< \alpha \leq 1,$ we give the following identities from \eqref{eq:def-sub}:
 \begin{equation}\label{eq:link-1}
{_{a}}I_t^{\alpha,\lambda}f(t)=e^{-\lambda t}{_{a}}I_t^{\alpha}[e^{\lambda t}f(t)],\qquad  {_{a}}D_t^{\alpha,\lambda}f(t)=e^{-\lambda t}{_{a}}D_t^{\alpha}[e^{\lambda t}f(t)].
 \end{equation}
Moreover, if $ f(a)=0,$ by the composition formula \cite{Podlubny1999}, it holds that
 \begin{equation}\label{link-2}
{_{a}}I_t^{\alpha,\lambda} {_{a}}D_t^{\alpha,\lambda}f(t)=f(t).
 \end{equation}

Based on the definition of the modified Gr\"{u}nwald derivative \cite{OldhamS1974} and the shifted Gr\"{u}nwald- Letnikov formula \cite{TadMeer2006}, we define
\begin{equation}\label{a1}
A_{\tau,r}^{\alpha,\lambda}f(t)=\frac{1}{\tau^{\alpha}}\sum_{k=0}^{+\infty} e^{-\lambda(k-r) \tau}g^{\alpha}_kf(t-(k-r)\tau),
\quad \forall~~ \alpha >0.
\end{equation}
where  $\tau$ is the step size, and $\{g^{\alpha}_k\}$  are the coefficients of the power series expansion of the function $(1-z)^{\alpha},$ i.e.,
\begin{equation}\label{a3}
(1-z)^{\alpha}=\sum_{k=0}^{+\infty}(-1)^k\binom{\alpha}{k}z^k=\sum_{k=0}^{+\infty} g^{\alpha}_kz^k,\quad \forall~~ |z|<1.
\end{equation}
Obviously, $g^{\alpha}_k= (-1)^k\binom{\alpha}{k}$ can also be evaluated recursively {\color{red} by}
\begin{equation}\label{recursion}
    g^{\alpha}_0=1,\quad g^{\alpha}_k=(1-\frac{\alpha+1}{k})g^{\alpha}_{k-1}, \quad k=1,2, \cdots.
\end{equation}

Let $f\in L_1{(\mathbb{R})}$.  Define the  function  space   {\color{red}$\mathscr{C}^{n+\alpha}(\mathbb{R})$ by }
$$\mathscr{C}^{n+\alpha}(\mathbb{R})=\left\{f| \int_{-\infty}^{+\infty}(1+|\omega|)^{n+\alpha}|\hat{f}(\omega)|d\omega<\infty \right\},$$
where $\hat{f}(\omega)=\int_{-\infty}^{+\infty}e^{i\omega t}f(t)dt$ is the  {\color{red}Fourier transform} of $f(t).$

Motivated by the work  in \cite{TadMeer2006}, we present an  asymptotic  expansion of $A_{\tau,r}^{\alpha,\lambda}f(t)$.

\begin{lem}\label{lem-asympototic-expansion}
Suppose that $f\in L_1{(\mathbb{R})}$ and   $f\in \mathscr{C}^{n+\alpha}(\mathbb{R}).$
Then
\begin{equation}\label{a2-1-5}
A_{\tau,r}^{\alpha,\lambda}f(t)={_{-\infty}}D_t^{\alpha,\lambda}f(t)+\sum_{k=1}^{n-1}c_{r,k}^{\alpha}\ {_{-\infty}}D_t^{\alpha+k,\lambda}f(t)\tau^k+\mathcal{O}(\tau^n)
\end{equation}
 uniformly  {\color{red} holds} in $t\in \mathbb{R}$ as $\tau\rightarrow 0,$ where $c_{r,k}^{\alpha}$ are the coefficients  of the power series expansion of the function ${\color{red}W_{\alpha,r}}(z)=(\frac{1-e^{-z}}{z})^{\alpha}e^{rz}.$
\end{lem}

\begin{proofed}
Using the {\color{red}Fourier transform}, we obtain
\begin{eqnarray*}
&&\mathcal{F}(A_{\tau,r}^{\alpha,\lambda}f )(\omega)=\frac{1}{\tau^{\alpha}}\sum_{k=0}^{\infty} e^{-\lambda(k-r) \tau}g^{\alpha}_k \mathcal{F}(f(t-(k-r)\tau))(\omega)\\
&&\qquad \qquad \qquad =\frac{1}{\tau^{\alpha}} \sum_{k=0}^{\infty} e^{-\lambda(k-r) \tau}g^{\alpha}_k e^{i\omega (k-r)\tau}\hat{f}(\omega)\\
&&\qquad \qquad \qquad =\frac{1}{\tau^{\alpha}}(1-e^{-(\lambda-i\omega) \tau} )^{\alpha} e^{r(\lambda-i\omega) \tau} \hat{f}(\omega)\\
&&\qquad \qquad \qquad =(\lambda-i\omega)^{\alpha} \biggl(  \frac{ 1-e^{-(\lambda-i\omega) \tau}}{(\lambda-i\omega)\tau }\biggl )^{\alpha} e^{r(\lambda-i\omega) \tau} \hat{f}(\omega)\\
&&\qquad \qquad \qquad =(\lambda-i\omega)^{\alpha} W_{\alpha,r}((\lambda-i\omega) \tau) \hat{f}(\omega).\\
\end{eqnarray*}
Observe that $c_{r,0}^\alpha=1$ and  $$\mathcal{F}({_{-\infty}}D_t^{\alpha+k,\lambda}f)(\omega)=(\lambda-i\omega)^{\alpha+k}\hat{f}(\omega) $$ for $k=0,1,2\cdots .$
Therefore, there exists
\begin{eqnarray}\label{sub-eq-2-1}
 \mathcal{F}(A_{\tau,r}^{\alpha,\lambda}f )(\omega)=\sum_{k=0}^{n-1}c_{r,k}^{\alpha}
 \mathcal{F}({_{-\infty}}D_t^{\alpha+k,\lambda}f)(\omega)\tau^k+\Phi(\omega,\tau),\quad
\end{eqnarray}
where $$\Phi(\omega,\tau)=(\lambda-i\omega)^{\alpha}[ W_{\alpha,r}((\lambda-i\omega) \tau)-\sum_{k=0}^{n-1}c_{r,k}^{\alpha}(\lambda-i\omega)^k \tau^k]\hat{f}(\omega).$$
Taking inverse Fourier {\color{red}transform} on the both sides of \eqref{sub-eq-2-1} leads to
\begin{equation}\label{sub-eq-2-2}
A_{\tau,r}^{\alpha,\lambda}f(t)-\sum_{k=0}^{n-1}c_{r,k}^{\alpha}\ {_{-\infty}}D_t^{\alpha+k,\lambda}f(t)\tau^k=\frac{1}{2\pi}\int_{-\infty}^{+\infty}\Phi(\omega,\tau)e^{-i\omega t}d\omega.
\end{equation}
Due to the fact  that $f\in \mathscr{C}^{n+\alpha}(\mathbb{R})$  and
${\color{red}W_{\alpha,r}}(z)=\sum_{k=0}^{n-1}c_{r,k}^{\alpha}z^k+O(z^n),$ then there exists a constant $c$ such that \begin{eqnarray*}
&& \quad A_{\tau,r}^{\alpha,\lambda}f(t)-{_{-\infty}}D_t^{\alpha,\lambda}f(t)-\sum_{k=1}^{n-1}c_{r,k}^{\alpha}\ {_{-\infty}}D_t^{\alpha+k,\lambda}f(t)\tau^k\nonumber\\
 &&\leq \frac{1}{2\pi}\int_{-\infty}^{+\infty}|\Phi(\omega,\tau)|d\omega\leq c\int_{-\infty}^{+\infty}(1+|\omega|)^{n+\alpha}|\hat{f}(\omega)|d\omega\cdot \tau^n=\mathcal{O}(\tau^n).
\end{eqnarray*}
for   sufficiently small $\tau.$ This ends the proof.
\end{proofed}

\begin{rem}
Take $\lambda=0,~\alpha=1$ in above lemma.  When  $r=0,$ \eqref{a2-1-5} reduces to the backward difference operator for the standard derivative:
$$f^{\prime}(t)=\frac{f(t)-f(t-\tau)}{\tau}+\mathcal{O}(\tau); $$
 when  $r=\frac{\alpha}{2}=\frac{1}{2},$  \eqref{a2-1-5} reduces to the central difference quotient operator:
 $$f^{\prime}(t)=\frac{f(t+\frac{1}{2}\tau)-f(t-\frac{1}{2}\tau)}{\tau}+\mathcal{O}(\tau^2). $$
 On the other hand, the accuracy can be improved by  shifting from the backward difference formula to the central difference counterpart.     From this point of  view, the fractional shift in the formula \eqref{a2-1-5}  is natural and reasonable to achieve high accuracy.
\end{rem}

From Lemma \ref{lem-asympototic-expansion}, we can take $A_{\tau,r}^{\alpha,\lambda}f(t)$ as a first-order approximation to the fractional substantial  derivative $ {_{-\infty}}D_t^{\alpha,\lambda}f(t),$ that is,
 \begin{equation*}
 A_{\tau,r}^{\alpha,\lambda}f(t)= {_{-\infty}}D_t^{\alpha,\lambda}f(t)+{\color{red}\mathcal{O}}(\tau),
 \end{equation*}
 which is consistent with the approximation of the Riemann-Liouville  derivative by the shifted Gr\"{u}nwald-Letnikov formula \cite{TadMeer2006}.
 Note that  $c_{r,1}^{
\alpha}=r-\frac{\alpha}{2}$.
We take $t=s-r\tau $ and  $r=\frac{\alpha}{2}$  in \eqref{a2-1-5}, then it follows that
\begin{equation*}
A_{\tau,r}^{\alpha,\lambda}f(s-\frac{\alpha}{2}\tau)={_{-\infty}}D_t^{\alpha,\lambda}f(s-\frac{\alpha}{2}\tau)+\mathcal{O}(\tau^2).
\end{equation*}
  By \eqref{a1}, and using the Lagrange linear interpolation at the points $t=s$ and $t=s-\tau $ for the approximation of  $ {_{-\infty}}D_t^{\alpha,\lambda}f(s-r\tau)$, and replacing $s$ with $t$,
we can derive the following second-order Gr\"{u}nwald-Letnikov-formula-based approximation for the fractional substantial derivative.
\begin{thm}\label{lem-6}
Suppose that $f\in L_1{(\mathbb{R})}$ and   $f\in \mathscr{C}^{2+\alpha}(\mathbb{R}).$
Then
\begin{equation*}
(1-\frac{\alpha}{2}) {_{-\infty}}D_t^{\alpha,\lambda}f(t)+\frac{\alpha}{2} {_{-\infty}}D_t^{\alpha,\lambda}f(t-\tau)=\frac{1}{\tau^{\alpha}}\sum_{k=0}^{+\infty}g_k^{\alpha,\lambda} f(t-k\tau)+\mathcal{O}(\tau^2)
\end{equation*}
{\color{red}uniformly holds} in $t\in \mathbb{R}$ as $\tau\rightarrow 0,$ where
$g_k^{\alpha,\lambda}=e^{-(k-\frac{\alpha}{2})\lambda \tau}g^{\alpha}_k.$
\end{thm}

To apply Theorem \ref{lem-6} to the generic function $f(t)$ $(f(0)= 0)$ on $[0,+\infty)$, we make the zero extension as follows:
\begin{equation*}
\tilde{f}(t)= \left \{\begin{array}{cc}
               f(t) ,& t\in  [0,+\infty),\\
               0,~&  \mbox{others}
              \end{array}
              \right.
\end{equation*}
 and further assume that  $\tilde{f}(t)$ {\color{red}satisfies} the assumptions in  Theorem \ref{lem-6}, i.e., $\tilde{f}(t)\in \mathscr{C}^{2+\alpha}(\mathbb{R}).$
Then it follows that
\begin{eqnarray}\label{formula-1}
\quad (1-\frac{\alpha}{2})\,{_{0}}D_t^{\alpha,\lambda}f(t)+\frac{\alpha}{2}\,{_{0}} D_t^{\alpha,\lambda}f(t-\tau) =\frac{1}{\tau^{\alpha}}\sum_{k=0}^{[\frac{t}{\tau}]}g^{\alpha,\lambda}_{k}f(t-k\tau)+\mathcal{O}(\tau^2).\qquad
\end{eqnarray}


\section{Derivation of the numerical scheme}
In this section, we  construct a numerical scheme for the problem \eqref{c1}-\eqref{cc1}. We first derive a  semi-discretized scheme by the proposed  approximation   \eqref{formula-1}. Subsequently,  with the use of  the fourth-order compact finite difference  approximation of the space  derivatives \cite{sun2005}, we give a fully discrete scheme. In what follows, we suppose that $u(\mathbf{x},0)=0$ without loss of  generality.   Otherwise, consider the equation with the solution   $v(\mathbf{x},t)=u(\mathbf{x},t)-e^{-\lambda t}u(\mathbf{x},0)$ instead.
 Furthermore, assume that the extended $u(\cdot,t)\in \mathscr{C}^{2+\alpha}(\mathbb{R}).$


\subsection{A semi-discrete  scheme}

Now we are in the position to discrete \eqref{c1} in time. Take the uniform time step size $\tau=T/N$ and let $t_n=n\tau,\;0\leq n\leq N.$ The temporal  domain $[0,T]$ is covered by $\Omega_{\tau}=\{t_n~|~0\leq n\leq N \}.$  For a given grid function $v=\{v^n\}_{n=0}^{N}$ defined  on $\Omega_{\tau},$ we  introduce   the average  operator:
\begin{equation}\label{denotion-2}
\mathcal{A}_t^{\alpha}v^n=(1-\frac{\alpha}{2})v^{n}+\frac{\alpha}{2}v^{n-1}
\end{equation}
and the discrete difference quotient  operator:
\begin{equation}\label{denotion-1}
 \delta_{t}^{\alpha,\lambda}v^n=\frac{1}{\tau^{\alpha}}\sum_{k=0}^{n}g^{\alpha,\lambda}_{k}v^{n-k}.
\end{equation}

 Considering the equation \eqref{c1} at the  grid points $(\mathbf{x},t_n),$ we have
\begin{equation}\label{c2}
{\color{red}{_0}D_t^{\alpha,\lambda} }u(\mathbf{x},t_n)=\Delta u(\mathbf{x},t_n)+F(\mathbf{x},t_n),\quad  0\leq n\leq N.
\end{equation}
Performing the average  operator $\mathcal{A}_t^{\alpha}$ defined in \eqref{denotion-2} on  the equation \eqref{c2},  we have
\begin{equation*}
\mathcal{A}_t^{\alpha}{\color{red}{_0}D_t^{\alpha,\lambda} }u(\mathbf{x},t_n)=\mathcal{A}_t^{\alpha}\Delta u(\mathbf{x},t_n)+\mathcal{A}_t^{\alpha}F(\mathbf{x},t_n),\quad 1\leq n\leq N.
\end{equation*}
Applying the formula  \eqref{formula-1} and  combining  \eqref{denotion-1}  yield
\begin{equation}\label{semi-discretize-1}
\delta_t^{\alpha,\lambda}u(\mathbf{x},t_n)=\mathcal{A}_t^{\alpha}\Delta u(\mathbf{x},t_n)+\mathcal{A}_t^{\alpha}F(\mathbf{x},t_n)+\mathcal{O}(\tau^2),\quad 1\leq n\leq N.
\end{equation}

 Noting that \eqref{semi-discretize-1} is a semi-discretization to  \eqref{c1},  in next subsection, we will derive a fully discretization by  the finite difference method. It is worth to point out  that \eqref{semi-discretize-1}  can also be solved by many other methods that we mentioned in the introduction part, such as spectral method, finite element method, {\color{red}discontinuous Galerkin  method,} etc.
\subsection{A fully discrete scheme}
In this part, we use the fourth-order compact {\color{red}approximation } for the spatial discretization of the problem \eqref{c1}-\eqref{cc1}. To fix the idea, we consider only one-dimensional problem \eqref{c1}-\eqref{cc1}  in this subsection and illustrate  the  {\color{red}generalization to} two-dimensional case via an example  in  Section 6.

Let $\Omega=(a,b)$. Take  a positive  integer $M$ and let
$h=(b-a)/M$,  $x_i=ih$, $0\leq i\leq M.$ Then the spatial domain $[a,b]$ is covered by $\Omega_{h}=\{x_i~|~0\leq i\leq M \}.$

For completeness of the paper, we introduce some results that are helpful for understanding the construction of the scheme.
The following lemma is required to achieve high-order accuracy  for spatial discretization.
\begin{lem}\label{spatial-1}\cite{sun2005}
Denote $\theta(s)=(1-s)^3[5- 3(1 - s)^2].$  If $G(x)\in  C^6[x_{i-1},x_{i+1}],$ $x_{i+1}=x_i+h,$ $x_{i-1}=x_i-h,$ it holds that
\begin{eqnarray*}
&&\frac{1}{12}[G^{\prime\prime}(x_{i-1})+10G^{\prime\prime}(x_i)+G^{\prime\prime}(x_{i+1})]=\frac{G(x_{i-1})-2
G(x_i)+G(x_{i+1})}{h^2}\\
&&\quad \quad \quad +\frac{h^4}{360}\int_0^1[G^{(6)}(x_i-sh)+G^{(6)}(x_i+sh)]\theta (s)ds,\quad 1\leq i\leq M-1.
\end{eqnarray*}
\end{lem}
Introduce the average  operator $\mathcal{A}_x$ in spatial direction  as
\begin{eqnarray*}
\mathcal{A}_xv_{i}=\left\{
\begin{array}{ll}
\frac{1}{12}(v_{i-1}+10v_{i}+v_{i+1}),~~~1\leq i\leq M-1,\\
v_{i},~~~~~~~~~~~~~~~~i=0~ \mbox{or}~ M.\\
\end{array}
\right.
\end{eqnarray*}
Denote by $\delta_{x}^2$  the second-order central difference quotient operator, that is  $$\delta_{x}^2v_i=\frac{1}{h^2}(v_{i-1}-2v_i+v_{i+1}).$$

Taking $\mathbf{x}=x_i$ and performing the average  operator $\mathcal{A}_x$ on both sides of  \eqref{semi-discretize-1} gives
\begin{eqnarray}\label{c4}
&&\mathcal{A}_x\delta_t^{\alpha,\lambda}u(x_i,t_n)=\mathcal{A}_x\mathcal{A}_t^{\alpha}\partial_{x}^{2}u(x_i,t_n)
+\mathcal{A}_x\mathcal{A}_t^{\alpha}{\color{red}F}(x_i,t_n)+\mathcal{O}(\tau^2),\nonumber \\
&&\quad 1\leq i\leq M-1,\quad 1\leq n\leq N.
\end{eqnarray}
Define the grid functions  $$U_i^n=u(x_i,t_n),\quad F^{n}_i=F(x_i,t_n),\quad   0\leq i\leq M,\quad 0\leq n\leq N. $$
Note that the operators $ \mathcal{A}_t^{\alpha}$ and $\delta_x^{2}$ can be commuted to each other. Then using Lemma \ref{spatial-1}  yields
\begin{eqnarray}
&&\mathcal{A}_x\delta_{t}^{\alpha,\lambda} U_i^{n}-\mathcal{A}_t^{\alpha}\delta_{x}^{2}U_i^{n}=\mathcal{A}_t^{\alpha}\mathcal{A}_x F_i^{n}+R_i^n,~ 1\leq i\leq M-1,~1\leq n\leq  N, \nonumber\\ \label{c55}
\end{eqnarray}
where there exists a constant $c_1$ such that
\begin{equation}\label{eq7}
|R_i^n|\leq c_1(\tau^2+h^4),\quad 1\leq i\leq M-1,\quad 1\leq n\leq N.
\end{equation}
Omitting the small terms $R_i^n$ in \eqref{c55},
denoting  $u_i^n$ as the approximation of $U_i^n$ and noticing the initial and boundary conditions
\begin{eqnarray}
&&U_i^0=0 ,\quad  0\leq i\leq M, \label{b5}   \\
&&U_0^n=\phi(x_0,{\color{red}t_n}),\quad U_M^n=\phi(x_M,{\color{red}t_n}),\quad 1\leq n\leq N,      \label{b6}
\end{eqnarray}
we obtain a $(\alpha,\lambda)-$dependent compact difference scheme
\begin{eqnarray}
&&\mathcal{A}_x \delta_{t}^{\alpha,\lambda} u_i^{n}-\mathcal{A}_t^{\alpha}\delta_{x}^{2}u_i^{n}=\mathcal{A}_t^{\alpha}\mathcal{A}_x F_i^{n},~ 1\leq i\leq M-1,~1\leq n\leq  N,\label{c5}\\
&&u_i^0=0 , ~  0\leq i\leq M  , \label{c19}   \\
&&u_0^n=\phi(x_0,{\color{red}t_n}),~ u_M^n=\phi(x_M,{\color{red}t_n}),\quad 1\leq n\leq N .     \label{c17}
\end{eqnarray}
\begin{rem}
When    $\lambda=\lambda(x), $  the problem \eqref{c1}-\eqref{cc1} becomes the Feynman-Kac equation considered in \cite{DengCB2015-JSC-substantial}.
In Section 5, we will show numerically that the adapted  scheme, i.e., taking $\lambda_i=\lambda(x_i)$ in above scheme,  is of second-order convergence in time and fourth-order convergence in space for the forward Feynman-Kac equation considered in \cite{DengCB2015-JSC-substantial}.
\end{rem}

\section{Stability and convergence analysis}
\setcounter{equation}{0}
In this section, we shall give the stability and convergence analysis for the scheme \eqref{c5}-\eqref{c17}.

Let
$V_h=\{v~|~v=(v_0,\cdots,v_M ) ~\mbox{is a grid function on}~ \Omega_h~\mbox{and } ~v_0=v_M=0\}.$
For any $u, v\in V_h,$ we denote the maximum norm  by $\displaystyle\|u\|_{\infty}=\max_{0\leq i\leq M }|u_i|$ and define the discrete inner products and induced norm as
$$(u,v)=h\sum_{i=1}^{M-1}u_iv_i, \quad  \|u\|=\sqrt{(u,u)}.  $$
Let
 \begin{equation}l_n^{\alpha}=\sum_{k=0}^n g_k^{\alpha},\quad n\geq 0. \label{SC-Notation-1} \end{equation}
Then $l_n^{\alpha}$  are the coefficients of the power series expansion of the function $(1-z)^{\alpha-1},$ i.e.,
\begin{equation*}
(1-z)^{\alpha-1}=\sum_{k=0}^{+\infty} l_{k}^{\alpha}z^k,
\end{equation*}
In fact, from $(1-z)^{\alpha}=(1-z)\cdot (1-z)^{\alpha-1},$ we have
\begin{equation*}\label{a7}
\sum_{k=0}^{+\infty} g^{\alpha}_{k}z^k=(1-z) \sum_{k=0}^{+\infty} l_{k}^{\alpha}z^k=l^{\alpha}_0+\sum_{k=1}^{+\infty} (l^{\alpha}_{k}-l^{\alpha}_{k-1} )z^k.
\end{equation*}
Comparing  the coefficients of power series expansion  on  both sides leads to the desired result.
In addition,  the coefficients $l_n^{\alpha}$ in \eqref{SC-Notation-1} satisfy the following properties for $0 < \alpha <1,$
\begin{equation}\label{aa5}
\left\{  \begin{array}{l}
           \displaystyle l_0^{\alpha}=1,\quad l_1^{\alpha}=1-\alpha   , \\
           \displaystyle
         l_0^{\alpha} \geq   l_1^{\alpha} \geq  l_2^{\alpha}\geq l_3^{\alpha}\geq \cdots \geq 0.
         \end{array}
 \right.
\end{equation}
And  it has the lower bound, i.e.,
\begin{equation}\label{SC-inequality-1}
l_{n-1}^{\alpha}\geq\frac{1 }{n^{\alpha} \Gamma{(1-\alpha) }}.
\end{equation}
 See Lemma 2.1 in \cite{DengCB2015-JSC-substantial} for details.


Prior  to  {\color{red}considering} the stability and  convergence of the scheme \eqref{c5}-\eqref{c17}, we need to provide  some auxiliary lemmas. Inspired by the work in \cite{Alikhanov2015}, we first give the following lemma.
\begin{lem}\label{lemma-1}
For a given  sequence $\{v^n\}_{n=0}^{\infty} ,$ we have
\begin{eqnarray}
&&v^{n}\delta_{t}^{\alpha,\lambda}v^n \geq \frac{1}{2} e^{ -\frac{\alpha}{2}\lambda\tau}\delta_{t}^{\alpha,2\lambda} (v^n)^2 +\frac{e^{- \frac{\alpha}{2}\lambda\tau}}{2l_0^{\alpha}} \tau^{\alpha}(\delta_{t}^{\alpha,\lambda}v^n)^2, \label{eq:2-4} \\
&& v^{n-1}\delta_{t}^{\alpha,\lambda}v^n \geq  \frac{1}{2} e^{(1- \frac{\alpha}{2})\lambda\tau} \delta_{t}^{\alpha,2\lambda} (v^n)^2-\frac{e^{(1- \frac{\alpha}{2})\lambda\tau}}{2(l_0^{\alpha}-l_1^{\alpha})}\tau^{\alpha}  (\delta_{t}^{\alpha,\lambda}v^n)^2\label{eq:2-5},
\end{eqnarray}
where  $$  \delta_{t}^{\alpha,2\lambda} (v^n)^2= \frac{1}{\tau^{\alpha}} \sum_{k=0}^{n} g_{k}^{\alpha,2\lambda} (v^{n-k})^2 . $$
\end{lem}
\begin{proofed}
Note that $g_k^{\alpha,\lambda}=e^{\lambda(\frac{\alpha}{2}-k)\tau}g_k^{\alpha}. $   The expressions  $\delta_{t}^{\alpha,\lambda}v^n,\quad \delta_{t}^{\alpha,2\lambda} (v^n)^2 $ can also be rewritten as
\begin{eqnarray}
&&\delta_{t}^{\alpha,\lambda}v^n=e^{\lambda (\frac{\alpha}{2}-n)\tau}\frac{1}{\tau^{\alpha}} \sum_{k=0}^nl_{n-k}^{\alpha} (e^{\lambda k\tau}v^{k}-e^{\lambda (k-1)\tau}v^{k-1} ) ,\label{eq:add42}\\
  &&\delta_{t}^{\alpha,2\lambda}(v^n)^2=e^{2\lambda (\frac{\alpha}{2}-n)\tau}\frac{1}{\tau^{\alpha}} \sum_{k=0}^nl_{n-k}^{\alpha} [(e^{\lambda k\tau}v^{k})^2-(e^{\lambda (k-1)\tau}v^{k-1})^2 ],\quad\label{eq:add43}
\end{eqnarray}
where $v^{-1}=0.$  To {\color{red} simplify} notation but without ambiguity, hereafter, we denote $v_{\lambda}^k:=e^{\lambda k\tau}v^{k}.$
Thus we have
\begin{eqnarray}
&& \quad e^{ (2n-\frac{\alpha}{2})\lambda\tau} [ v^{n}\delta_{t}^{\alpha,\lambda}v^n - \frac{1}{2}e^{ -\frac{\alpha}{2}\lambda\tau} \delta_{t}^{\alpha,2\lambda} (v^n)^2 ] \nonumber \\
&& =v_{\lambda}^n \frac{1}{\tau^{\alpha}} \sum_{k=0}^nl_{n-k}^{\alpha} (v_{\lambda}^{k}-v_{\lambda}^{k-1} )- \frac{1}{\tau^{\alpha}} \sum_{k=0}^nl_{n-k}^{\alpha} (v_{\lambda}^{k}-v_{\lambda}^{k-1} ) \left(\frac{v_{\lambda}^{k}+v_{\lambda}^{k-1}}{2} \right)\nonumber \\
&& = \frac{1}{\tau^{\alpha}} \sum_{k=0}^nl^{\alpha}_{n-k} (v_{\lambda}^{k}-v_{\lambda}^{k-1} ) \left(v_{\lambda}^n-\frac{v_{\lambda}^{k}+v_{\lambda}^{k-1}}{2} \right)\nonumber \\
&& = \frac{1}{\tau^{\alpha}} \sum_{k=0}^nl^{\alpha}_{n-k} (v_{\lambda}^{k}-v_{\lambda}^{k-1} ) \biggl[\frac{v_{\lambda}^{k}-v_{\lambda}^{k-1}}{2}+\sum_{j=k+1}^{n}(v_{\lambda}^{j}-v_{\lambda}^{j-1}) \biggl]\nonumber \\
&& = \frac{1}{2\tau^{\alpha}} \sum_{k=0}^nl^{\alpha}_{n-k} (v_{\lambda}^{k}-v_{\lambda}^{k-1} )^2+ \frac{1}{\tau^{\alpha}}\sum_{j=1}^{n} (v_{\lambda}^{j}-v_{\lambda}^{j-1})  \sum_{k=0}^{j-1}  l^{\alpha}_{n-k} (v_{\lambda}^{k}-v_{\lambda}^{k-1} ) .\nonumber \\ \label{eq:2-6}
\end{eqnarray}
Denote
\begin{equation}
 w_{\lambda}^{j}=  \sum_{k=0}^{j-1}  l^{\alpha}_{n-k} (v_{\lambda}^{k}-v_{\lambda}^{k-1} ),\quad 1\leq j\leq n \label{eq:2-8}
\end{equation}
and
\begin{equation}
 w_{\lambda}^{n+1}=\sum_{k=0}^{n}  l^{\alpha}_{n-k} (v_{\lambda}^{k}-v_{\lambda}^{k-1} )=\tau^{\alpha}e^{ (n-\frac{\alpha}{2})\lambda\tau}\delta_{t}^{\alpha,\lambda}v^n . \label{eq:2-9}
\end{equation}
Then it follows that
\begin{eqnarray}\label{eq:4-1}
v_{\lambda}^{0}-v_{\lambda}^{-1}= \frac{w_{\lambda}^{1}}{l^{\alpha}_n}; \qquad v_{\lambda}^{j}-v_{\lambda}^{j-1}= \frac{w_{\lambda}^{j+1}-w_{\lambda}^{j} }{l^{\alpha}_{n-j} } ,\quad 1\leq j\leq n.
\end{eqnarray}
Substituting \eqref{eq:2-8}-\eqref{eq:4-1} into  \eqref{eq:2-6} leads to
\begin{eqnarray}
&& \quad e^{ (2n-\frac{\alpha}{2})\lambda\tau} [ v^{n}\delta_{t}^{\alpha,\lambda}v^n - \frac{1}{2}e^{ -\frac{\alpha}{2}\lambda\tau} \delta_{t}^{\alpha,2\lambda} (v^n)^2 ] \nonumber \\
&&=\frac{1}{2\tau^{\alpha}}l_n^{\alpha} ( \frac{w_{\lambda}^{1}}{l^{\alpha}_n} )^2 +\frac{1}{2\tau^{\alpha}} \sum_{k=1}^n l^{\alpha}_{n-k} \biggl(\frac{w_{\lambda}^{k+1}-w_{\lambda}^{k} }{l^{\alpha}_{n-k} }\biggl)^2 +\frac{1}{\tau^{\alpha}}\sum_{j=1}^{n} \biggl(\frac{w_{\lambda}^{j+1}-w_{\lambda}^{j} }{l^{\alpha}_{n-j} }\biggl) w_{\lambda}^{j}\nonumber \\
&&=\frac{1}{2\tau^{\alpha}} \frac{(w_{\lambda}^{n+1})^2}{l^{\alpha}_0}+\frac{1}{2\tau^{\alpha}} \sum_{k=1}^{n}\biggl(\frac{1}{l^{\alpha}_{n+1-k}}-\frac{1}{l^{\alpha}_{n-k}} \biggl)(w_{\lambda}^{k})^2. \label{eq:2-10}
\end{eqnarray}
Multiplying  \eqref{eq:2-10} by $e^{ (\frac{\alpha}{2}-2n)\lambda\tau}$ and  recalling   \eqref{aa5}, we have
\begin{eqnarray*}
&& \quad v^{n}\delta_{t}^{\alpha,\lambda}v^n - \frac{1}{2}e^{ -\frac{\alpha}{2}\lambda\tau} \delta_{t}^{\alpha,2\lambda} (v^n)^2 \geq \frac{e^{ (\frac{\alpha}{2}-2n)\lambda\tau}}{2\tau^{\alpha}} \frac{(w_{\lambda}^{n+1})^2}{l^{\alpha}_0}= \frac{e^{ -\frac{\alpha}{2}\lambda\tau}}{2l^{\alpha}_0} \tau^{\alpha}(\delta_{t}^{\alpha,\lambda}v^n)^2.
\end{eqnarray*}
Next, we prove the inequality  \eqref{eq:2-5}. Considering  the difference,  using the notations \eqref{eq:2-8}-\eqref{eq:2-9} again and noting the equality \eqref{eq:2-10},  we have
\begin{eqnarray*}
&& \quad e^{ (2n-1-\frac{\alpha}{2})\lambda\tau} [v^{n-1}\delta_{t}^{\alpha,\lambda}v^n - \frac{1}{2} e^{ (1-\frac{\alpha}{2})\lambda\tau} \delta_{t}^{\alpha,2\lambda} (v^n)^2+\frac{e^{(1- \frac{\alpha}{2})\lambda\tau}}{2(l_0^{\alpha}-l_1^{\alpha})}\tau^{\alpha}  (\delta_{t}^{\alpha,\lambda}v^n)^2] \nonumber \\
&&=e^{ (2n-\frac{\alpha}{2})\lambda\tau} [v^{n}\delta_{t}^{\alpha,\lambda}v^n - \frac{1}{2}e^{ -\frac{\alpha}{2}\lambda\tau}  \delta_{t}^{\alpha,2\lambda} (v^n)^2 ] +\frac{\tau^{\alpha}}{2(l_0^{\alpha}-l_1^{\alpha})}  e^{ (2n-\alpha)\lambda\tau}(\delta_{t}^{\alpha,\lambda}v^n)^2 \nonumber \\
&&\quad ~~- [e^{n\lambda \tau}v^{n}-e^{(n-1)\lambda \tau}v^{n-1}]e^{ (n-\frac{\alpha}{2})\lambda\tau}\delta_{t}^{\alpha,\lambda}v^n \nonumber \\
&&=\frac{1}{2\tau^{\alpha}} \frac{(w_{\lambda}^{n+1})^2}{l_0^{\alpha}}+\frac{1}{2\tau^{\alpha}} \sum_{k=1}^{n}\biggl(\frac{1}{l_{n+1-k}^{\alpha}}-\frac{1}{l_{n-k}^{\alpha}} \biggl)(w_{\lambda}^{k})^2+ \frac{\tau^{\alpha}}{2(l_0^{\alpha}-l_1^{\alpha})} ( \frac{w_{\lambda}^{n+1}}{\tau^{\alpha}})^2\nonumber \\
&&\quad ~~- \frac{(w_{\lambda}^{n+1}-w_{\lambda}^{n})w_{\lambda}^{n+1} }{l^{\alpha}_{0} \tau^{\alpha}} \nonumber \\
&&= \frac{1}{\tau^{\alpha}}  \biggl( \frac{1}{2(l_0^{\alpha}-l_1^{\alpha})}-\frac{1}{2l_0^{\alpha}} \biggl)(w_{\lambda}^{n+1})^2+ \frac{1}{\tau^{\alpha}}\frac{1}{l_{0}^{\alpha}}w_{\lambda}^nw_{\lambda}^{n+1}+\frac{1}{2\tau^{\alpha}}\biggl(\frac{1}{l_{1}^{\alpha}}-\frac{1}{l^{\alpha}_{0}} \biggl)(w_{\lambda}^{n})^2 \nonumber\\
&&\qquad  +\frac{1}{2\tau^{\alpha}} \sum_{k=1}^{n-1}\biggl(\frac{1}{l_{n+1-k}^{\alpha}}-\frac{1}{l_{n-k}^{\alpha}} \biggl)(w_{\lambda}^{k})^2\nonumber 
\end{eqnarray*}
\begin{eqnarray*}
&&= \frac{1}{2l_0^{\alpha} \tau^{\alpha}}\biggl( \sqrt{\frac{l_1^{\alpha}}{l_0^{\alpha}-l_1^{\alpha}}}w_{\lambda}^{n+1} +\sqrt{\frac{l_0^{\alpha}-l_1^{\alpha}}{l_1^{\alpha}}}w_{\lambda}^{n}  \biggl)^2  +\frac{1}{2\tau^{\alpha}} \sum_{k=1}^{n-1}\biggl(\frac{1}{l^{\alpha}_{n+1-k}}-\frac{1}{l^{\alpha}_{n-k}} \biggl)(w_{\lambda}^{k})^2\nonumber \\
&& \geq 0,
\end{eqnarray*}
which completes the proof.
 \end{proofed}
\begin{lem}\label{SC-lemma-1}
Given a sequence $\{v^n\} $ we have
\begin{eqnarray}
\mathcal{A}_t^{\alpha} v^{n}\cdot  \delta_{t}^{\alpha,\lambda}v^n &\geq& \frac{1}{2} {\color{red}\left[\frac{\alpha}{2}e^{\lambda \tau}+(1-\frac{\alpha}{2})\right] } e^{-\frac{\alpha}{2}\lambda \tau} \delta_{t}^{\alpha,2\lambda} (v^n)^2\nonumber\\
&&+\frac{2-\alpha-e^{\lambda \tau}}{4}\tau^{\alpha}e^{-\frac{\alpha}{2}\lambda \tau}(\delta_{t}^{\alpha,\lambda}v^n)^2. \label{eq:2-7}
\end{eqnarray}
In particular, when $\alpha=1,~\lambda=0$ it holds that
$$\mathcal{A}_t v^{n}\cdot\delta_{t}v^n: =\frac{v^n+v^{n-1}}{2} \cdot  \frac{v^n-v^{n-1}}{\tau} = \frac{(v^n)^2-(v^{n-1})^{2}}{2\tau} . $$
\end{lem}
\begin{proofed}
From  Lemma \ref{lemma-1}, we have
\begin{eqnarray*}
&& \mathcal{A}_t^{\alpha} v^{n}\cdot  \delta_{t}^{\alpha,\lambda}v^n =(1-\frac{\alpha}{2}) v^n\delta_{t}^{\alpha,\lambda}v^n+\frac{\alpha}{2}v^{n-1} \delta_{t}^{\alpha,\lambda}v^n\nonumber \\
&& \geq (1-\frac{\alpha}{2})e^{-\frac{\alpha}{2}\lambda \tau}\left[\frac{1}{2} \delta_{t}^{\alpha,2\lambda} (v^n)^2 +\frac{\tau^{\alpha}}{2l_0^{\alpha}} (\delta_{t}^{\alpha,\lambda}v^n)^2\right]\nonumber \\
&&\qquad +\frac{\alpha}{2}e^{(1-\frac{\alpha}{2})\lambda \tau}\left[\frac{1}{2} \delta_{t}^{\alpha,2\lambda} (v^n)^2-\frac{\tau^{\alpha}}{2(l_0^{\alpha}-l_1^{\alpha})} (\delta_{t}^{\alpha,\lambda}v^n)^2\right]\nonumber \\
&&= \frac{1}{2}{\color{red}\left[\frac{\alpha}{2}e^{\lambda \tau}+(1-\frac{\alpha}{2})\right] } e^{-\frac{\alpha}{2}\lambda \tau}\delta_{t}^{\alpha,2\lambda} (v^n)^2 \nonumber \\ &&\qquad +\left[(1-\frac{\alpha}{2})\frac{1}{2l_0^{\alpha}}-\frac{\alpha}{2}\frac{e^{\lambda \tau}}{2(l^{\alpha}_0-l^{\alpha}_1)}\right]\tau^{\alpha}e^{-\frac{\alpha}{2}\lambda \tau}(\delta_{t}^{\alpha,\lambda}v^n)^2.
\end{eqnarray*}
Substituting  $l^{\alpha}_0=1$ and $ l^{\alpha}_1=1-\alpha  $ into above inequality  leads to \eqref{eq:2-7}, 
which completes the proof.
\end{proofed}
\begin{lem}\cite{Sunbook2009}\label{lem4}
For any $u,v\in V_h,$ it holds that
$(\mathcal{A}_x u,v)=(u,\mathcal{A}_x v).$
\end{lem}
\begin{lem}\cite{Sunbook2009}\label{SC-lemma-5}
For any $v\in V_h,$ it holds that
$$\frac{2}{3}\|v\|^2 \leq (\mathcal{A}_x v,v)\leq \|v\|^2.$$
\end{lem}
Since $\mathcal{A}_x$ is positive definite and
self-adjoint, we can consider its square {\color{red}root} denoted as $Q_x$, and have
\begin{equation}\label{eq:add41}
(\mathcal{A}_x u,v)=(Q_x u,Q_xv).
\end{equation}
Naturally, we define an equivalent  norm  as
\begin{equation*}
 \|v\|_A= \sqrt{(\mathcal{A}_xv,v)}=\sqrt{(Q_xv,Q_xv)}.
\end{equation*}

\begin{lem}\label{lem3}\cite{Sunbook2009}
For any $v\in V_h,$ it holds that
$-(\delta_{x}^2v, v)\geq C_{\Omega}\|v\|^2,$ where $C_\Omega$ is a constant independent to the step size $h$ but related to the domain $\Omega$.
\end{lem}

Now we give a prior estimate for the difference scheme \eqref{c5}-\eqref{c17}.
\begin{thm}\label{SC-thm1}$\mathbf{(Priori~ estimate)}$
Suppose  $\{v_i^n\}$ be the solution of
\begin{eqnarray}
&&\mathcal{A}_x\delta_{t}^{\alpha,\lambda} v_i^{n}-\mathcal{A}_t^{\alpha}\delta_{x}^2v_i^{n}=S_i^{n},~ 1\leq i\leq M-1,~1\leq n\leq  N,\label{c10} \\
&&v^n_0=0,\quad v^n_M=0,\quad 1 \leq n\leq N,  \label{c11}  \\
&&v^0_i=v_0(x_i),\quad 0\leq i\leq M. \label{c12}
\end{eqnarray}
If  $\;2-\alpha-e^{\lambda \tau}\geq 0,$ then
\begin{equation}
\|v^n\|^2\leq \frac{3}{2}e^{2|\lambda|T}\|v^{0}\|^2+\frac{3\Gamma{(1-\alpha)}T^{\alpha}e^{2|\lambda|T}}{C_{\Omega}} \max_{0\leq n\leq N}\|S^n\|^2 \label{SC-inequality-6}
\end{equation}
for $1<n\leq N,$ where $v_0(x_0)=v_0(x_M)=0,\quad \displaystyle \|S^n\|=\sqrt{h\sum_{i=1}^{M-1}(S_i^n)^2 }. $
\end{thm}
\begin{proofed}
Taking the inner product of \eqref{c10} with  $\mathcal{A}_t^{\alpha}v^{n},$  we have
\begin{equation}\label{c13}
 (\mathcal{A}_x\delta_{t}^{\alpha,\lambda}v^{n},\mathcal{A}_t^{\alpha}v^{n})-(\mathcal{A}_t^{\alpha}\delta_{x}^2v^{n},\mathcal{A}_t^{\alpha}v^{n})=(S^n,\mathcal{A}_t^{\alpha}v^{n}).
\end{equation}
By \eqref{eq:add41}, for the first term on the left-hand side of \eqref{c13}, we get
\begin{equation*}
    (\mathcal{A}_x\delta_{t}^{\alpha,\lambda}v^{n},\mathcal{A}_t^{\alpha}v^{n})=(Q_x\delta_{t}^{\alpha,\lambda}v^{n},Q_x\mathcal{A}_t^{\alpha}v^{n})  =(\delta_{t}^{\alpha,\lambda}Q_xv^{n},\mathcal{A}_t^{\alpha}Q_xv^{n})    .
\end{equation*}
By Lemma \ref{SC-lemma-1} and {\color{red} assumption $\;2-\alpha-e^{\lambda \tau}\geq 0$,} we have
\begin{eqnarray}\label{c14}
   \quad (\delta_{t}^{\alpha,\lambda}Q_xv^{n},\mathcal{A}_t^{\alpha}Q_xv^{n})  &=& h\sum_{i=1}^{M-1}\delta_{t}^{\alpha,\lambda}(Q_xv^{n})_i\cdot  \mathcal{A}_t^{\alpha}(Q_xv^{n})_i\nonumber\\
   &\geq&  \frac{1}{2}e^{-\frac{\alpha}{2}\lambda \tau} {\color{red}\left[\frac{\alpha}{2}e^{\lambda \tau}+(1-\frac{\alpha}{2})\right] } h\sum_{i=1}^{M-1}\delta_{t}^{\alpha,2\lambda} (Q_xv^{n})_i^2.\nonumber\\
\end{eqnarray}
For the second term on the left-hand side,
from Lemma \ref{lem3}, we obtain
\begin{equation}\label{c15}
-(\mathcal{A}_t^{\alpha}\delta_{x}^2v^{n},\mathcal{A}_t^{\alpha}v^{n})=-(\delta_{x}^2\mathcal{A}_t^{\alpha}v^{n},\mathcal{A}_t^{\alpha}v^{n})\geq C_{\Omega} \|\mathcal{A}_t^{\alpha}v^{n}\|^2.
\end{equation}
As to the term on the right-hand side, in view of Lemma \ref{lem4} we have
\begin{equation}\label{c16}
(S^n,\mathcal{A}_t^{\alpha}v^{n})\leq \|S^n\|\cdot \|\mathcal{A}_t^{\alpha}v^{n}\|\leq \frac{1}{4C_{\Omega}}\|S^n\|^2+C_{\Omega} \|\mathcal{A}_t^{\alpha}v^{n}\|^2.
\end{equation}
Substituting \eqref{c14}-\eqref{c16} into \eqref{c13}, we obtain
\begin{equation*}
  \frac{1}{2}e^{-\frac{\alpha}{2}\lambda \tau}{\color{red}\left[\frac{\alpha}{2}e^{\lambda \tau}+(1-\frac{\alpha}{2})\right] } h\sum_{i=1}^{M-1}\delta_{t}^{\alpha,2\lambda} (Q_xv^{n})_i^2 \leq \frac{1}{4C_{\Omega}}\|S^n\|^2.
\end{equation*}
{\color{red}Since $\displaystyle\frac{\alpha}{2}e^{\lambda \tau}+(1-\frac{\alpha}{2})\geq \frac{1}{2} ,$ we have
\begin{equation*}
 e^{-\frac{\alpha}{2}\lambda \tau}h\sum_{i=1}^{M-1}\delta_{t}^{\alpha,2\lambda} (Q_xv^{n})_i^2 \leq \frac{1}{C_{\Omega}}\|S^n\|^2.
\end{equation*}
}
Consequently, by \eqref{eq:add43}, it follows that
\begin{eqnarray}
&&\quad l_0^{\alpha}\|v_{\lambda}^{n}\|_A^2 \nonumber \\
&&\leq \sum_{k=1}^{n-1}(l^{\alpha}_{n-k-1}-l^{\alpha}_{n-k} )  \|v_{\lambda}^{k}\|_A^2 +   (l^{\alpha}_{n-1}-l^{\alpha}_n) \|v_{\lambda}^{0}\|_A^2 + \frac{\tau^{\alpha}}{C_{\Omega}} e^{(2n-\frac{\alpha}{2})\lambda \tau}\|S^n\|^2
\nonumber \\
&&\leq  \sum_{k=1}^{n-1}(l^{\alpha}_{n-k-1}-l^{\alpha}_{n-k} )  \|v_{\lambda}^{k}\|_A^2+ l^{\alpha}_{n-1}(\|v_{\lambda}^{0}\|_A^2+\frac{\tau^{\alpha}}{C_{\Omega}l^{\alpha}_{n-1}} e^{(2n-\frac{\alpha}{2})\lambda \tau}\|S^n\|^2)
\nonumber \\
&&\leq  \sum_{k=1}^{n-1}(l^{\alpha}_{n-k-1}-l^{\alpha}_{n-k} )  \|v_{\lambda}^{k}\|_A^2+ l^{\alpha}_{n-1}(\|v_{\lambda}^{0}\|^2+\frac{\Gamma{(1-\alpha)}T^{\alpha}}{C_{\Omega}} e^{(2n-\frac{\alpha}{2})\lambda \tau}\|S^n\|^2),\nonumber \\
&& \qquad\qquad 1\leq n\leq N,  \label{SC-inequality-4}
\end{eqnarray}
where we have used  the estimate \eqref{SC-inequality-1}  in the last inequality.
Denote $$E^N=\|v_{\lambda}^{0}\|^2+{\color{red}\frac{\Gamma{(1-\alpha)}T^{\alpha}}{C_{\Omega}}\max_{1\leq n\leq N}e^{(2n-\frac{\alpha}{2})\lambda \tau} \max_{1\leq n\leq N}\|S^n\|^2} .$$
The inequality  \eqref{SC-inequality-4} is simplified to
\begin{eqnarray}
&&l_0^{\alpha}\|v_{\lambda}^{n}\|_A^2\leq    \sum_{k=1}^{n-1}(l^{\alpha}_{n-k-1}-l^{\alpha}_{n-k} )  \|v_{\lambda}^{k}\|_A^2+ l^{\alpha}_{n-1}E^N,
\quad 1\leq n\leq N. \nonumber\\\label{SC-inequality-5}
\end{eqnarray}
Next we prove the following inequality
\begin{eqnarray}
&&\|v_{\lambda}^{n}\|_A^2\leq   E^N,
\quad 1\leq n\leq N. \label{sub-4-1}
\end{eqnarray} by the mathematical induction method.  Obviously, it follows   from  \eqref{SC-inequality-5} that the equality \eqref{sub-4-1} holds for $n=1.$  Suppose  that \eqref{sub-4-1} still holds for $n=1,2,\cdots m.$
From \eqref{SC-inequality-5} at $n=m+1$, we have
\begin{eqnarray*}
&&l_0^{\alpha}\|v_{\lambda}^{m+1}\|_A^2\leq    \sum_{k=1}^{m}(l^{\alpha}_{m-k}-l^{\alpha}_{m+1-k} )  \|v_{\lambda}^{k}\|_A^2+ l^{\alpha}_{m}E^N \nonumber \\
&& \qquad \qquad \quad   \leq    \sum_{k=1}^{m}(l^{\alpha}_{m-k}-l^{\alpha}_{m+1-k} )  E^N+ l^{\alpha}_{m}E^N =l^{\alpha}_0E^N,
\end{eqnarray*}
which leads to \eqref{sub-4-1}.

It follows from Lemma \ref{SC-lemma-5}  that  $$\frac{2}{3}e^{2\lambda n\tau}\|v^{n}\|^2=\frac{2}{3}\|v_{\lambda}^{n}\|^2\leq \|v_{\lambda}^{n}\|_A^2\leq E^N. $$
Thus
\begin{eqnarray*}
&&\|v^n\|^2 \leq \frac{3}{2} e^{-2\lambda n\tau}E^N \\
&\leq &\frac{3}{2} e^{-2\lambda n\tau}\biggl( \|v^{0}\|^2+\frac{\Gamma{(1-\alpha)}T^{\alpha}}{C_{\Omega}}\max_{1\leq n\leq N}e^{(2n-\frac{\alpha}{2})\lambda \tau} \max_{1\leq n\leq N}\|S^n\|^2 \biggl)\\
&\leq &\frac{3}{2} e^{2|\lambda|T} \|v^{0}\|^2+\frac{3\Gamma{(1-\alpha)}T^{\alpha}}{C_{\Omega}}e^{2|\lambda|T } \max_{1\leq n\leq N}\|S^n\|^2,\;1\leq n\leq N.
\end{eqnarray*}

This completes the proof.
\end{proofed}

Using Theorem \ref{SC-thm1}, it is readily to obtain the following result.
\begin{thm}$\mathbf{(Stability)}$
The difference scheme \eqref{c5}-\eqref{c17} is unconditionally stable {\color{red} with respect to } the initial value  and right-hand term.
\end{thm}

We now turn to the convergence of the difference scheme \eqref{c5}-\eqref{c17}.

\begin{thm}$\mathbf{(Convergence)}$  Assume that  {\color{red}$u(x, t) \in C^{6,3}(\Omega \times [0, T]) $}  be the solution of the problem \eqref{c1}-\eqref{cc1},
 and $\{u_i^n\}$ be the solution of difference scheme \eqref{c5}-\eqref{c17}. We further assume that zero extended function  $u(x,t) $ satisfies  the assumptions in Theorem \ref{lem-6}.  Let $e_i^n=u(x_i,t_n)-u_i^n,\quad 0\leq i\leq M,~0\leq n\leq N. $
Then   the following estimate
$$\|e^n\|\leq \sqrt{\frac{3\Gamma{(1-\alpha)e^{2|\lambda|T}}T^{\alpha}(b-a)}{C_{\Omega}}}c_1(\tau^2+h^4),\quad 1\leq n\leq N $$
holds, where  $C_\Omega$ is defined as that in Lemma \ref{lem3}.
\end{thm}
\begin{proofed}
Subtracting \eqref{c5}-\eqref{c17} from  \eqref{c55}, \eqref{b5}-\eqref{b6} respectively, we get the error equations
\begin{eqnarray*}
&&\mathcal{A}_x\delta_{t}^{\alpha,\lambda}e_i^{n}-\mathcal{A}_t^{\alpha}\delta_{x}^2e_i^{n}=R_i^n, \quad 1\leq i\leq M-1,\quad 1\leq n\leq N,\\
&&e^n_0=0,\quad e^n_M=0,\quad 1 \leq n\leq N,   \\
&&e^0_i=0,\quad 0\leq i\leq M.
\end{eqnarray*}
It follows from Theorem \ref{SC-thm1} and  the truncation  error estimate \eqref{eq7} that
\begin{eqnarray*}
\|e^n\|^2 &\leq& \frac{3\Gamma{(1-\alpha)}e^{2|\lambda|T}T^{\alpha}}{C_{\Omega}} \max_{0\leq n\leq N}\|R^n\|^2 \\
&\leq& \frac{3\Gamma{(1-\alpha)}e^{2|\lambda|T}T^{\alpha}(b-a)}{C_{\Omega}}c_1^2(\tau^2+h^4)^2, \quad 1\leq n\leq N.
\end{eqnarray*}
This completes the proof.
\end{proofed}

\section{An improved algorithm}\label{sec5}


This section presents an improved algorithm for the problem \eqref{c1}-\eqref{cc1} with nonsmooth solution.
To begin with this section, we first introduce the following  lemma.
\begin{lem}\label{sub-lem-5-2}
Let $f(t)=e^{-\lambda t}t^{\beta},\,\,\beta >0.$
Then
\begin{equation*}
(1-\frac{\alpha}{2}) [{_{-\infty}}D_t^{\alpha,\lambda}f(t)]_{t=t_n}+\frac{\alpha}{2} [{_{-\infty}}D_t^{\alpha,\lambda}f(t)]_{t=t_{n-1}}=\frac{1}{\tau^{\alpha}}\sum_{k=0}^{n}g_k^{\alpha,\lambda} f(t_{n-k})+\mathcal{O}(t_n^{\beta-\alpha-2}\tau^2)
\end{equation*}
uniformly {\color{red} holds} in $t\in \mathbb{R}$ as $\tau\rightarrow 0,$ where
$g_k^{\alpha,\lambda}=e^{-(k-\frac{\alpha}{2})\lambda \tau}g^{\alpha}_k.$
\end{lem}
The proof of above lemma is similar to \cite{Lubich1986} {\color{red}and} here we skip it.

To fix the idea, let us first consider
 the fractional ordinary differential equation with substantial derivative as below
\begin{equation}\label{sub-eq-5-3}
{_0}D_t^{\alpha,\lambda}u(t)=\mu u(t)+F(t).
\end{equation}
We can seek its solution of the following form
\begin{equation}\label{sub-eq-5-1}
u(t)=\sum_{j,k=0}^{+\infty}c_{j,k}e^{-\lambda t}t^{k+j\alpha},
\end{equation}
which can be derived by the method of  power series.
Indeed, suppose that the right-hand side function $F$ is sufficiently smooth and  has the power series expansion of the form below
\begin{equation}\label{sub-eq-5-2}
F(t)=\sum_{k=0}^{+\infty}F_{k}e^{-\lambda t}t^{k}.
\end{equation}
Inserting \eqref{sub-eq-5-1} and \eqref{sub-eq-5-2} into \eqref{sub-eq-5-3} leads to
\begin{eqnarray*}
\sum_{j,k=0}^{+\infty}c_{j,k}{_0}D_t^{\alpha,\lambda}(e^{-\lambda t}t^{k+j\alpha})=\mu\sum_{j,k=0}^{+\infty}c_{j,k}e^{-\lambda t}t^{k+j\alpha}+\sum_{k=0}^{+\infty}F_{k}e^{-\lambda t}t^{k}.
\end{eqnarray*}
Consequently
\begin{eqnarray}\label{sub-eq-5-4}
&&\sum_{j,k=0}^{+\infty}c_{j,k}e^{-\lambda t} \frac{\Gamma(k+1+j\alpha)}{\Gamma(k+1+(j-1)\alpha)}t^{k+(j-1)\alpha}\nonumber\\
&&=\mu\sum_{j,k=0}^{+\infty}c_{j,k}e^{-\lambda t}t^{k+j\alpha}+\sum_{k=0}^{+\infty}F_{k}e^{-\lambda t}t^{k},
\end{eqnarray}
where we have used the formula
$$ {_0}D_t^{\alpha,\lambda}(e^{-\lambda t}t^{k+j\alpha})=e^{-\lambda t} \frac{\Gamma(k+1+j\alpha)}{\Gamma(k+1+(j-1)\alpha)}t^{k+(j-1)\alpha}.$$
Comparing the coefficients of both sides of identity \eqref{sub-eq-5-4} can solve  $c_{j,k},$ the coefficients to be determined.

Now we turn our attention back to the TFSDE \eqref{c1}-\eqref{cc1}.
From the aforementioned  discussion,    it makes sense  to  suppose  the solution has the following  structure:
\begin{equation}\label{SC-5-1}
u(\mathbf{x},t)=\sum_{j,k=0}^{+\infty}c_{j,k}(\mathbf{x})e^{-\lambda(\mathbf{x}) t}t^{k+j\alpha}.
\end{equation}
From the asymptotic behavior of the solution near the original point $t=0$ in \eqref{SC-5-1}, we can clearly see, by Lemma \ref{sub-lem-5-2},  that the expected convergence order of the   scheme proposed in the previous section may decrease correspondingly and thus the accuracy will deteriorate meantime. To make up for the {\color{red}lost} accuracy, we have to turn to other approaches.
One of the natural ideas occurs to us: if we can  remove or separate  the singular part    from the candidate solution, then the remaining part has desired regularity to satisfy the necessary smoothness assumptions.  However, this seems to be ideal but  can not work in practice since the coefficients $c_{j,k}$ is unknown and difficult to calculate. An efficient approach firstly proposed by Lubich \cite{Lubich1986}   is to add  starting quadrature to make sure that the corrected scheme is accurate for the power function $t^{\beta},$  $\beta \in \Theta,$ where $\Theta=\{\beta~|~ \beta=k+j\alpha: j,\,k\in \mathbb{N}, \, \beta \leq 2+\alpha  \} .$
 Following the spirit of Lubich \cite{Lubich1986}, by Lemma \ref{sub-lem-5-2}  we get the following  semi-discrete  correction scheme
  \begin{equation}\label{semi-discretize-correct}
\delta_t^{\alpha,\lambda}u(\mathbf{x},t_n)+\sum_{k=0}^{S}w_{n,k}^{\alpha,\lambda}u(\mathbf{x},t_k)=\mathcal{A}_t^{\alpha}\Delta u(\mathbf{x},t_n)+\mathcal{A}_t^{\alpha}F(\mathbf{x},t_n)+\mathcal{O}(\tau^2),\quad 1\leq n\leq N,
\end{equation}
for the problem \eqref{c1}-\eqref{cc1} with solution in the form of \eqref{SC-5-1},  where $S$ is the number of elements of the  set $\Theta$ and  the  starting weights $ w_{n,k}^{\alpha,\lambda}$ can be computed from the Vandermonde type system
\begin{eqnarray}
&&\sum_{k=0}^{S}w_{n,k}^{\alpha,\lambda}t_k^{\beta_j}e^{-\lambda t_k}=(1-\frac{\alpha}{2})\frac{\Gamma(\beta_j+1)}{\Gamma(\beta_j+1-\alpha)}t_n^{\beta_j-\alpha}e^{-\lambda t_n}\nonumber\\
&& +\frac{\alpha}{2} \frac{\Gamma(\beta_j+1)}{\Gamma(\beta_j+1-\alpha)} t_{n-1}^{\beta_j-\alpha}e^{-\lambda t_{n-1}} -\frac{1}{\tau^{\alpha}}\sum_{k=0}^{n}g_k^{\alpha,\lambda} t_n^{\beta_j}e^{-\lambda t_k},\;\,j=0,\cdots, S.\quad\qquad\label{eq:57}
\end{eqnarray}
Instead of studying the theory of the improved algorithm, we will provide a two-dimensional example to show its validity in next  section, where a fully discretized compact finite  scheme for the two-dimensional TFSDE is also given; see Example 6.3.

{\color{red}\begin{rem}
Based on the  pitfalls in fast numerical solvers
for fractional differential equations studied in \cite{Diethelm2006}, it is known that the linear system \eqref{eq:57} is ill-posed, which may result in the reduced efficiency of the algorithm \eqref{semi-discretize-correct}. However, using a few correction terms can significantly increase the accuracy for problems with low regularity. It will not be computational expensive to calculate the necessary starting weights; see also \cite{Zengnew}.  
\end{rem}}

\section{Numerical examples}

 In this section, we give some numerical results to illustrate the performance of the proposed scheme and verify  the theoretical prediction. In Example 6.1, a prototypical equation is considered  to show the convergence order of the proposed second-order approximation \eqref{formula-1}  with or without the smoothness assumptions in Theorem \ref{lem-6}.   Example 6.2 provides a comparison between the proposed scheme \eqref{c5}-\eqref{c17} and the scheme in \cite{DengCB2015-JSC-substantial}, for solving the backward Feynman-Fac equation \eqref{eq:feynman-kac}. In Example 6.3, we illustrate efficiency of the improved algorithm by solving a two-dimensional TFSDE with nonsmooth solution.

\begin{example}\label{example-3}
Consider the following prototypical  equation
\begin{eqnarray*}
&&_0D_t^{\alpha,\lambda} u(t)=f(t),
\end{eqnarray*}
with $u(0)$ and  the source term $f(t) $ such that  the equation admits the solution
$u(t)=e^{-\lambda t}(t^3+t^\nu ) ,$  and  $ 0<\alpha<1.$ We take $\lambda=0.5$ and change $\nu$ to test the convergence order of the approximation \eqref{formula-1} for the problem with different regularity.
\end{example}
Denote
$$E(\tau)=\max_{1\leq n \leq N}|u^n-U^n|.$$ {\color{red}In this example  $Rate: =\log_2{\frac{E(\tau)}{E(\tau/2)}}.$}

\begin{table}[!thb]
\centering
\caption{\color{red} Error and convergence order of the approximation \eqref{formula-1} for $u(t)$ with different regularity (Example \ref{example-3}).  }\label{table-1}
\vspace{0.1in}
\scalebox{0.85}{\begin{tabular}{c|c|cc|cc|cc}
\hline
 & &  \multicolumn{2}{c|}{ $\alpha=0.2$}  & \multicolumn{2}{c|}{ $\alpha=0.5$} &\multicolumn{2}{c}{ $\alpha=0.8$}  \\
\cline{3-4} \cline{5-6} \cline{7-8}
&$\tau$& $E(\tau)$  & Rate &$E(\tau)$  & Rate &$E(\tau)$  & Rate \\
\hline
$\nu=2.5$ & 1/16    & 1.9225e-4 &      & 4.7884e-4 &      & 7.5901e-4 & \\
          & 1/32    & 4.8101e-5 & 2.00 & 1.2002e-4 & 1.99 & 1.9083e-4 & 2.00\\
          & 1/64    & 1.2029e-5 & 2.00 & 3.0043e-5 & 2.00 & 4.7871e-5 & 2.00\\
          & 1/128   & 3.0076e-6 & 2.00 & 7.5152e-6 & 2.00 & 1.1993e-5 & 2.00\\
\hline
$\nu=2$ & 1/16    & 1.5725e-4 &      & 3.8400e-4 &      & 5.6264e-4 & \\
        & 1/32    & 3.9392e-5 & 2.00 & 9.6827e-5 & 1.99 & 1.4297e-4 & 1.98\\
        & 1/64    & 9.8586e-6 & 2.00 & 2.4348e-5 & 1.99 & 3.6235e-5 & 1.98\\
        & 1/128   & 2.4661e-6 & 2.00 & 6.1116e-6 & 1.99 & 9.1649e-6 & 1.98\\
\hline
$\nu=1.5$ & 1/16    & 1.3578e-4 &      & 3.1355e-4 &      & 3.0475e-4 & \\
          & 1/32    & 3.6556e-5 & 1.89 & 7.8536e-5 & 2.00 & 1.2018e-4 & 1.34\\
          & 1/64    & 1.2691e-5 & 1.53 & 1.9651e-5 & 2.00 & 4.4176e-5 & 1.44\\
          & 1/128   & 4.4625e-6 & 1.51 & 4.9148e-6 & 2.00 & 1.5849e-5 & 1.48\\
\hline
$\nu=1$ & 1/16    & 8.0608e-4 &      & 3.1355e-4 &      & 1.0800e-3 & \\
        & 1/32    & 4.0503e-4 & 0.99 & 7.0492e-4 & 1.00 & 5.2210e-4 & 1.05\\
        & 1/64    & 2.0355e-4 & 0.99 & 3.5386e-4 & 0.99 & 2.6071e-4 & 1.00\\
        & 1/128   & 1.0211e-4 & 1.00 & 1.7745e-4 & 1.00 & 1.3082e-4 & 0.99\\
\hline
$\nu=0.5$ & 1/16    & 1.0492e-2 &      & 2.7597e-2 &      & 4.9525e-2 & \\
          & 1/32    & 7.5289e-3 & 0.48 & 1.9804e-2 & 0.48 & 3.5543e-2 & 0.48\\
          & 1/64    & 5.3646e-3 & 0.49 & 1.4111e-2 & 0.49 & 2.5327e-2 & 0.49\\
          & 1/128   & 3.8081e-3 & 0.49 & 1.0017e-2 & 0.49 & 1.7978e-2 & 0.49\\
\hline
\end{tabular}}
\end{table}

Tables \ref{table-1} displays   that the approximation \eqref{formula-1} can obtain second-order accuracy in uniform maximum norm when $\nu=2.5$ and $2$, while  the convergence order decays with the decrease of $\nu$ evidently. {\color{red}The reduction of accuracy implies that some additional assumptions are needed to get second-order convergence for the approximation \eqref{formula-1}, e.g., the vanishing derivatives $u_t(0)=0$.  To keep the second-order convergence of the  approximation \eqref{formula-1}, we make an assumption of  the extended solution $\tilde{u}(t)\in \mathscr{C}^{2+\alpha}(\mathbb{R})$ in this paper. To some extent the assumption may be too strict, so we have considered the use of correction terms in the proposed schemes for problems with low regularity; see Section \ref{sec5} for the improved schemes and Example \ref{example-2D} for the corresponding numerical results.  }


  \begin{example}\label{example-2}
Consider the backward fractional Feynman-Kac equation presented in \cite{DengCB2015-JSC-substantial}
\begin{eqnarray}\label{eq:feynman-kac}
{_0}\partial_t^{\alpha,\lambda(x)}[P(x,\rho,t)-e^{-\rho U(x)t}P(x,\rho,0) ] =k_{\alpha} \partial_x^2P(x,\rho,t)+f(x,t) \quad
\end{eqnarray}
on a finite domain $0<x<1,$ $ 0<t\leq 1,$ with the coefficients $\kappa_{\alpha}=0.5$, and $\lambda(x)=\rho U(x),$ $U(x)=x, $ $\rho=1+\sqrt{-1} ;$
the forcing function
\begin{eqnarray}
 &&f(x,t)=-\kappa_{\alpha} e^{-\rho x t}(t^{3+\alpha}+1)(\rho^2 t^2 \sin(\pi x )-2\pi\rho t\cos(\pi x)-\pi^2\sin(\pi x) ), \nonumber\\
&& \qquad \qquad +\frac{\Gamma(4+\alpha )}{\Gamma(4)} e^{ -\rho xt}t^3\sin(\pi x).
  \end{eqnarray}
the initial condition $P(x,\rho,0)=\sin (\pi x), $ and the boundary conditions ${\color{red}P}(0,\rho, t)={\color{red}P}(1,\rho,t)=0. $
The exact solution is given by
$ P(x,\rho,t)= e^{-\rho x t}(t^{3+\alpha} +1)\sin (\pi x). $
\end{example}
Notice that the initial value is not equal to zero. So we have made the transformation
\begin{equation*}
u(x,t)=P(x,\rho,t)-e^{-\rho xt}P(x,\rho,0).
\end{equation*}

\begin{table}[!thb]
\centering
\caption{ Error and convergence order of the scheme \eqref{c5}-\eqref{c17} in time with a fixed  space step-size $h=1/40$  and $\lambda(x)=(1+\sqrt{-1})x$ (Example \ref{example-2}). }\label{table5}
\vspace{2mm}
\newsavebox{\tableboxwan}
\begin{lrbox}{\tableboxwan}
\begin{tabular}{c|cc|cc|cc}
\hline
 & \multicolumn{2}{c|}{ $\alpha=0.2$}  & \multicolumn{2}{c|}{ $\alpha=0.5$} &\multicolumn{2}{c}{ $\alpha=0.8$}  \\
\cline{2-3} \cline{4-5} \cline{6-7}
   $\tau$ & $E_1(\tau,h) $ & $Rate$ & $E_1(\tau,h) $ & $Rate$ &$E_1(\tau,h) $ & $Rate$  \\
     \hline
    1/5           & 5.1150e-3     &              & 1.3602e-2     &            & 2.1793e-2     &          \\
    1/10          & 1.3025e-3     &     1.97     & 3.4574e-3     &     1.98             & 5.5020e-3     &     1.99    \\
     1/20            & 3.2845e-4     &     1.99         & 8.7149e-4     &     1.99       & 1.3819e-3     &     1.99 \\
     1/40          & 8.2375e-5     &     2.00        & 2.1870e-4     &     1.99         & 3.4619e-4     &     2.00        \\
\hline
\end{tabular}
\end{lrbox}
\scalebox{0.85}{\usebox{\tableboxwan}} 
\end{table}

\begin{table}[!thb]
\centering
\caption{ Error and convergence order of the scheme \eqref{c5}-\eqref{c17} in space  with a fixed   time step $\tau=1/10000$  and $\lambda(x)=(1+\sqrt{-1})x $ (Example \ref{example-2}). }\label{table6}
\vspace{2mm}
\newsavebox{\tableboxrong}
\begin{lrbox}{\tableboxrong}
\begin{tabular}{c|cc|cc|cc}
\hline
 & \multicolumn{2}{c|}{ $\alpha=0.2$}  & \multicolumn{2}{c|}{ $\alpha=0.5$} &\multicolumn{2}{c}{ $\alpha=0.8$}  \\
\cline{2-3} \cline{4-5} \cline{6-7}
   $h$ & $E_1(\tau,h) $ & $Rate$ & $E_1(\tau,h) $ & $Rate$ &$E_1(\tau,h) $ & $Rate$  \\
     \hline
    1/4          & 2.0681e-03     &              & 1.7956e-03     &              & 1.4439e-03     &           \\
    1/8         & 1.2566e-04     &     4.04  & 1.1059e-04     &     4.02              & 9.0999e-05     &     3.99    \\
     1/16             & 7.9126e-06     &     3.99          & 6.9117e-06     &     4.00        & 5.6411e-06     &     4.01\\
     1/32          & 4.9271e-07     &     4.01      & 4.2904e-07     &     4.01         & 3.4901e-07     &     4.01       \\
     1/64 & 2.9927e-08     &     4.04  & 2.4756e-08     &     4.12      & 1.8798e-08     &     4.21 \\
\hline
\end{tabular}
\end{lrbox}
\scalebox{0.85}{\usebox{\tableboxrong}} 
\end{table}
Actually,  replacing the average operators  $\mathcal{A}_t^{\alpha}$ and  $\mathcal{A}_x$ in \eqref{c5}-\eqref{c19} by identity operators in temporal  and spatial  directions   respectively, we can obtain the following  scheme in \cite{DengCB2015-JSC-substantial} as below
\begin{eqnarray}
&&\delta_{t}^{\alpha,\lambda_i} u_i^{n}-\delta_{x}^{2}u_i^{n}=F_i^{n},~ 1\leq i\leq M-1,~1\leq n\leq  N,\label{Deng-3c5}\\
&&u_i^0=u_0(x_i) , ~  0\leq i\leq M, \label{Deng-3c19}   \\
&&u_0^n=\phi_a(x_0,{\color{red}{t_n}}),~ u_M^n=\phi_b(x_M,{\color{red}{t_n}}),\quad 1\leq n\leq N.     \label{Deng-3c17}
\end{eqnarray}

\begin{figure}
  \centerline{\includegraphics[width=8cm]{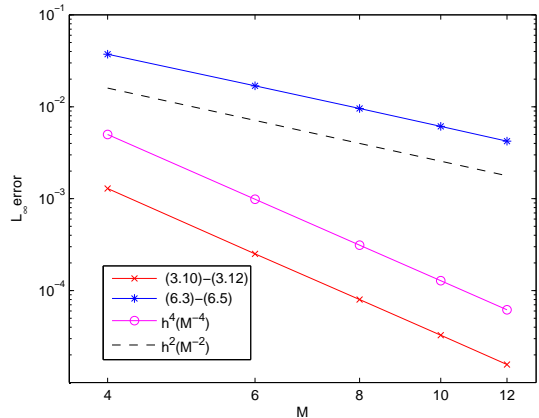}}
\caption{ Convergence order and error comparison between      the scheme \eqref{c5}-\eqref{c17}  and the scheme \eqref{Deng-3c5}-\eqref{Deng-3c17} in \cite{DengCB2015-JSC-substantial} with $\alpha=0.5$ for solving Equation \eqref{eq:feynman-kac} at $t=1$. In all cases, we set $N=M^2$  (Example \ref{example-2}).}
\label{fig2}
\end{figure}

Define the error $$E_1(\tau,h)=\max_{1\leq i\leq M-1}|u(x_i,t_N)-u_{i}^N|,$$
where $u(x_i,t_N) $ represents the exact solution and $u_{i}^N$ is the  numerical solution with the mesh step-sizes $h $ and $\tau$ at the grid point $ (x_i, t_N).$
Assume
$$E_1(\tau,h) = O(\tau ^p)+O(h^q).$$ If $\tau$ is sufficiently small, then $E_1(\tau,h) \approx O(h^q)$. Consequently,
$\frac{E_1(\tau,2h)}{E_1(\tau,h)}\approx 2^q$ and  $q \approx  \log_2\left(\frac{E_1(\tau,2h)}{E_1(\tau,h)}\right)$ is the convergence order with respect to the spatial step size. {\color{red}Likewise, the convergence order in time can be taken as  $p \approx  \log_2\left(\frac{E_1(\tau,2h)}{E_1(\tau,h)}\right)$  with   sufficiently small $h.$}
Tables 6.2-6.3 list the errors and convergence orders of our scheme,  showing that it is of second-order convergence   in time and fourth-order convergence   in space respectively. This is in a good agreement with our theoretical results.
Evidently, we can see from Figure \ref{fig2} that our proposed  scheme displays  much higher accuracy and enjoys  advantage over the scheme \eqref{Deng-3c5}-\eqref{Deng-3c19} presented  in \cite{DengCB2015-JSC-substantial}.
Table \ref{table7} further   shows that when the two schemes generate the same accuracy, our proposed scheme needs fewer temporal and spatial grid points and less CPU time than that of scheme \eqref{Deng-3c5}-\eqref{Deng-3c19}.

\begin{table}[!htb]
\centering
\caption{A comparison between the  scheme \eqref{c5}-\eqref{c17}  and the scheme \eqref{Deng-3c5}-\eqref{Deng-3c17} in \cite{DengCB2015-JSC-substantial} for different $\alpha$ and step sizes for solving Equation \eqref{eq:feynman-kac} (Example \ref{example-2}).}\label{table7}
\vspace{2mm}
\newsavebox{\tableboxzh}
\begin{lrbox}{\tableboxzh}
\begin{tabular}{c|cccc|cccc}
\hline
 & \multicolumn{4}{c|}{the scheme \eqref{c5}-\eqref{c17}}  & \multicolumn{4}{c}{ the scheme \eqref{Deng-3c5}-\eqref{Deng-3c17} }   \\
\cline{2-5} \cline{6-9}
   &$N$ & $M $ & $E_1(\tau,h) $ & CPU(second) &$N$ & $M$ &$E_1(\tau,h) $ & CPU(second)  \\
     \hline
             &16       &4        & 1.9804e-3     &   0.0151             &    256   &  16     & 2.2831e-3     &   1.0427      \\
       $\alpha=0.1$      &36       &6      & 3.8109e-4     &   0.0340              &   1296   & 36    & 4.5106e-4     &  50.9968     \\
          &64       &8        & 1.1946e-4     &   0.0661                &   4096   & 64         & 1.4269e-4     & 1484.3471\\
\hline
             &16       &4        & 1.2906e-3     &   0.0090             &    256   &  16     & 2.3822e-3     &   0.9892    \\
$\alpha=0.5$      &36       &6      & 2.5103e-4     &   0.0262             &   1296   & 36    & 4.6994e-4     &  50.6735     \\
           &64       &8        & 7.9988e-5     &   0.0671                &   4096   & 64        & 1.4863e-4     & 1491.8938  \\
             \hline
             &16       &4       & 1.6037e-3     &   0.0077             &    256   &  16      & 2.8498e-3     &   0.9932      \\
     $\alpha=0.9$        &36       &6      & 3.4057e-4     &   0.0261              &   1296   & 36    & 5.6214e-4     &  51.8453    \\
 &64       &8       & 1.0771e-4     &   0.0687                &   4096   & 64       & 1.7779e-4     & 1471.3060 \\
             \hline
\end{tabular}
\end{lrbox}
\scalebox{0.85}{\usebox{\tableboxzh}} 
\end{table}

\begin{example}\label{example-2D}
Consider the following two dimensional  backward fractional\\ Feynman-Kac equations
\begin{eqnarray}\label{eq:feynman-kac2D}
{_0}\partial_t^{\alpha,\lambda(x,y)}[u(x,y,t)-e^{-\lambda(x,y)t}u(x,y,0) ] = \Delta u(x,y,t)+F(x,y,t) \quad
\end{eqnarray}
on a finite domain $\Omega=(0,1)\times(0,1),$ $ 0<t\leq 1,$  and $\lambda(x,y)=0.01(x+y).$
We take the initial condition $u(x,y,0)=\sin(\pi x)\sin(\pi y), $ and the boundary condition $u(x,y, t)=0$ for $(x,y)\in\partial \Omega$.
The exact solution is
\begin{eqnarray*}
 &&u(x,y,t)= e^{ -\lambda(x,y)t}(1+t^{\alpha}+t^{2\alpha}+t^3)\sin(\pi x)\sin(\pi y).
  \end{eqnarray*}
\end{example}

\begin{table}[!t]\label{table-c-1}
\centering
\caption{ Error and  convergence order of the scheme \eqref{2c5}-\eqref{2c17} for the 2D problem \eqref{eq:feynman-kac2D}. $\alpha=0.2$, $h_1=h_2=1/60$; $S$ denotes the number of correction terms (Example \ref{example-2D}).}
\vspace{2mm}
\newsavebox{\tableboxro}
\begin{lrbox}{\tableboxro}
\begin{tabular}{c|cc|cc|cc|cc}
\hline
 & \multicolumn{2}{c|}{ $ S=0$}  &\multicolumn{2}{c|}{ $S=1$}&\multicolumn{2}{c|}{ $S=2$} &\multicolumn{2}{c}{ $S=3$} \\
\cline{2-3} \cline{4-5} \cline{6-7}\cline{8-9}
   $N$ & $E_2(\tau,h_1,h_2) $ & $Rate$ & $E_2(\tau,h_1,h_2) $ & $Rate$ &$E_2(\tau,h_1,h_2) $ & $Rate$ &$E_2(\tau,h_1,h_2) $ & $Rate$ \\
     \hline
 4       & 9.32e-3     &              & 8.86e-4     &           & 6.68e-4     &       & 1.47e-3     &   \\
          8& 8.68e-3     &     0.10      & 2.76e-4     &     1.68  & 1.95e-4     &     1.78 & 2.34e-4     &     2.65 \\
 16      & 8.03e-3     &     0.11       & 2.02e-4     &     0.45     & 5.09e-5     &     1.94 & 5.31e-5     &     2.14\\
 32   & 7.43e-3     &     0.11    & 1.70e-4     &     0.25    & 1.32e-5     &     1.95 & 1.33e-5     &     2.00 \\
\hline
\end{tabular}
\end{lrbox}
\scalebox{0.8}{\usebox{\tableboxro}} 
\end{table}
For  spatial approximation, take  two integers $M_1,M_2 $ and let
$h_1=1/M_1,$ $h_2=1/M_2,$
   $x_i=ih_1,$ $0\leq i\leq M_1,$ $y_j=j h_2, $ $0\leq j\leq M_2.$
Let $\bar{\Omega}_h=\{(x_i,y_j)|0\leq i\leq M_1,~ 0\leq j\leq M_2\},$ and $\Omega_h=\bar{\Omega}_h\cap\Omega, $ and $\partial \Omega_h=\bar{\Omega}_h\cap\partial\Omega. $
For any grid function $v=\{v_{i,j}|0\leq i\leq M_1, 0\leq j\leq M_2\},$ denote
$$\delta_xv_{i-\frac{1}{2},j}=\frac{1}{h_1}(v_{i,j}-v_{i-1,j} ),\quad \delta_x^2 v_{i,j}=\frac{1}{h_1}(\delta_xv_{i+\frac{1}{2},j}-\delta_xv_{i-\frac{1}{2},j}). $$
Similar notations $\delta_yv_{i,j-\frac{1}{2}}, $ $\delta_y^2v_{i,j} $ can be defined.
The spatial average  operators are defined as
\begin{eqnarray*}
\mathcal{A}_xv_{i,j}=\left\{
\begin{array}{ll}
\frac{1}{12}(v_{i-1,j}+10v_{i,j}+v_{i+1,j}),~~~1\leq i\leq M_1-1,~~0\leq j\leq M_2,\\
v_{i,j},~~~~~~~~~~~~~~~~i=0~ \mbox{or}~ M_1,~~0\leq j\leq M_2,\\
\end{array}
\right.\\
\mathcal{A}_yv_{i,j}=\left\{
\begin{array}{ll}
\frac{1}{12}(v_{i,j-1}+10v_{i,j}+v_{i,j+1}),~~~1\leq j\leq  M_2-1,~~0\leq i\leq M_1,\\
v_{i,j},~~~~~~~~~~~~~~~~j=0 ~\mbox{or} ~M_2,~~0\leq i\leq M_1.\\
\end{array}
\right.
\end{eqnarray*}
By generalizing the improved scheme \eqref{semi-discretize-correct}, we present a compact difference  scheme with correction terms in time for the two-dimensional equation  \eqref{eq:feynman-kac2D}
\begin{eqnarray}
&&\mathcal{A}_x \mathcal{A}_y[\delta_{t}^{\alpha,\lambda} u_{i,j}^{n}+\sum_{k=0}^{S}w_{n,k}^{\alpha,\lambda}u_{i,j}^{n}]-\mathcal{A}_t^{\alpha}\mathcal{A}_y\delta_{x}^{2}u_{i,j}^{n}-\mathcal{A}_t^{\alpha}\mathcal{A}_x\delta_{y}^{2}u_{i,j}^{n}= \mathcal{A}_t^{\alpha}\mathcal{A}_x\mathcal{A}_y F_{i,j}^{n},\nonumber \\
&&\quad (x_i,y_j)\in \Omega_h, \quad ~1\leq n\leq  N,\label{2c5}\\
&&u_{i,j}^0=u_0(x_i,y_j) , ~\quad    (x_i,y_j)\in \bar{\Omega}_h, \label{2c19}   \\
&&u_{i,j}^n=\phi(x_i,y_j,{\color{red}{t_n}}),\quad (x_i,y_j)\in \partial \Omega_h,\quad 1\leq n\leq N .     \label{2c17}
\end{eqnarray}
Denote
$$E_2(\tau,h_1,h_2)=\max_{1\leq n\leq N}\|{\color{red}u(\cdot,\cdot,t_n)-u^n} \|_{\infty}.$$

Table 6.5 shows that the errors  behave well and convergence order is augmented  up to two  with the increase of the number of correction terms, which indicates the efficiency and validity of our improved algorithm.

\section{Conclusion}
We have  proposed a second-order   approximation for  the fractional substantial derivative  and applied it to solve the time-fractional substantial diffusion equation. By combining the fourth-order compact finite difference approximation, a fully discrete Gr\"{u}nwald-Letnikov-formula-based compact difference scheme   has been  presented.  It has been  proved that the proposed scheme is unconditionally stable and  of second-order convergence in time and fourth-order convergence in space.

 We  have numerically showed that our proposed scheme can reach the predicted accuracy when solving a problem with smooth solution. Compared  with the scheme proposed in \cite{DengCB2015-JSC-substantial} for solving the backward Feymann-Fac fractional equation,   our scheme  requires  less storage and computational cost for the same accuracy.  While if the solution is not smooth enough, we have illustrated in Example \ref{example-3} that our proposed scheme loses its accuracy to some extent. To overcome this difficulty, we  followed the idea of Lubich \cite{Lubich1986} and introduced the corresponding  starting quadratures of substantial version. We have proposed an improved algorithm by applying the correction terms to increase  the accuracy of the scheme near the origin and recover the second-order accuracy in time for problems with nonsmooth solution,   which has been demonstrated by Example 6.3.

By the modified Gr\"{u}nwald derivative and the shifted Gr\"{u}nwald-Letnikov-formula, we have presented a weighted and shifted substantial Gr\"{u}nwald formula and given its asymptotic expansion, which is essential to derive the second-order approximation for the fractional substantial derivative. The proposed weighted average method is an efficient tool to get high-order scheme, not limited to the  problem under consideration in this study. {\color{red} It  can also be extended and applied to solve tempered fractional equations, which have attracted considerable interest recently; see \cite{ChenDeng2014-arxiv-substantial-tempered,Zayernouri2015}.} It should be also pointed out that  our scheme is consistent with the Crank-Nicolson scheme when the  fractional derivatives  reduces to the  standard operators, which is one of  attractive features of our proposed scheme.

%

\section*{Acknowledgement}
The research   is partially
supported by National Natural Science Foundation of China (No.
11271068) and by the Fundamental Research Funds for the Central Universities and the Research and Innovation Project for College Graduates of Jiangsu Province (Grant No.: KYLX\_0081).
G. Lin would like to thank the support by NSF Grant DMS-1115887, and the U.S. Department of Energy, Office of Science, Office of Advanced Scientific Computing Research, Applied Mathematics program as part of the Collaboratory on Mathematics for Mesoscopic Modeling of Materials, and Multifaceted Mathematics for Complex Energy Systems.

\end{document}